%% file: nonproper-nse-loewner.tex
\documentclass[%
  twocolumn
   , colorlinks 
]{mpi2015-cscpreprint}
\usepackage{cite}
\usepackage{amsmath,amssymb,amsfonts}
\usepackage{algorithmic}
\usepackage{graphicx}
\usepackage{textcomp}
\usepackage{xcolor}
\usepackage[]{hyperref}

\input{my_macros}

\def\BibTeX{{\rm B\kern-.05em{\sc i\kern-.025em b}\kern-.08em
    T\kern-.1667em\lower.7ex\hbox{E}\kern-.125emX}}

\definecolor{matlabblue}{HTML}{0072BD}
\definecolor{matlaborange}{HTML}{D95319}
\definecolor{matlabyellow}{HTML}{EDB120}
\definecolor{matlabpurple}{HTML}{7E2F8E}
\definecolor{matlabgreen}{HTML}{77AC30}
\definecolor{matlablightblue}{HTML}{4DBEEE}
\definecolor{matlabred}{HTML}{A2142F}

\def\tt{^\intercal}
\def\bHpoly{\mathbf{H}_{\mathsf{poly}}}
\def\bHspr{\mathbf{H}_{\mathsf{spr}}}

\hypersetup{
  colorlinks=true,
  urlcolor=matlabblue,
  citecolor=matlabgreen,
  linkcolor=matlabred
}

\begin{document}


\title{Implicit and explicit matching of non-proper transfer functions in the Loewner framework
}


\author[$\ast$]{Ion Victor Gosea}
\affil[$\ast$]{Max Planck Institute for Dynamics of Complex Technical Systems,
  Magdeburg, Germany. 
  \authorcr
  Email: \texttt{\href{mailto:gosea@mpi-magdeburg.mpg.de}{gosea@mpi-magdeburg.mpg.de}},
  ORCID: \texttt{\href{https://orcid.org/0000-0003-0228-8522}%
    {0000-0003-3580-4116}}}

\author[$\ast\ast$]{Jan Heiland}
\affil[$\ast\ast$]{Max Planck Institute for Dynamics of Complex Technical Systems,
  Magdeburg, Germany. \newline
Faculty of Mathematics, Otto von Guericke University Magdeburg, Germany.
  \authorcr
  Email: \texttt{\href{mailto:heiland@mpi-magdeburg.mpg.de}{heiland@mpi-magdeburg.mpg.de}},
  ORCID: \texttt{\href{https://orcid.org/0000-0003-0228-8522}%
    {0000-0003-0228-8522}}}

\shortauthor{I.V. Gosea and J. Heiland}
\shortdate{}
  

  \abstract{
The reduced-order modeling of a system from data (also known as system identification) is a classical task in system and control theory and well understood for standard linear systems with the
so-called Loewner framework as one of many established approaches. In the case of descriptor systems for which the transfer function is not proper anymore, recent
research efforts have addressed strategies to deal with the non-proper parts more or less
explicitly. In this work, we propose a variant of a Loewner matrix-based interpolation
algorithm that implicitly addresses possibly non-proper components of the system response. We evaluate the performance of the suggested approach by comparing against recently-developed explicit algorithms for which we propose a linearized Navier-Stokes model with a significant non-proper behavior as a benchmark
example.
}

\novelty{An extension of the \emph{Antoulas-Anderson} algorithm to implicitly
identify non-proper parts of the transfer function and a flow simulation
benchmark example with a significant linear term.}

\maketitle

\section{Introduction}

Accurate modeling of physical phenomena often leads to large-scale dynamical
systems that require long simulation times and storage of large amount of data.
In this context, model order reduction (MOR) aims at obtaining much smaller and
simpler models that are still capable of accurately representing the behavior of
the original process. The Loewner framework (LF)
\cite{antoulas2017tutorial}
is very appealing due to its data-driven nature, which makes it non-intrusive as it
does not use the full or exact description of the model. Hence, it can be viewed as a data-driven reduced-order modeling tool (in this context, data are frequency response measurements). 

Dynamical systems characterized by differential algebraic equations (DAEs) (referred to as descriptor, singular or semi-state systems) are not only of theoretical interest but also have a broad application range. In chemistry, for example, the additional algebraic equations account for thermodynamic equilibrium relations, steady state assumptions, or empirical correlations \cite{pantelides1988mathematical}. In mechanics, DAEs result from holonomic and non-holonomic constraints \cite{morMehS05}. We refer the reader to \cite{kunkel2006differential} for an in-depth account of analysis and numerical solution of DAEs, and to contributions from the last two decades that extended classical MOR methods to specific cases of dynamical systems with DAEs in \cite{morSty04,morMehS05,morReiS10,morBenS17}.

The LF is based on the
Loewner pencil that allows solving the generalized realization problem for
linear time-invariant (LTI) systems \cite{mayo2007framework}, and obtaining
reduced-order models through compression based on the Singular Value Decomposition (SVD). The Loewner matrix method of Antoulas-Anderson (AA) in \cite{AA86} also uses a Loewner matrix to construct interpolatory rational functions, but it is based on barycentric representations of such functions, as opposed to the Loewner pencil formulation in LF.
 To broaden the applicability of LF, extensions of the LF have recently been proposed for specific classes of dynamical systems characterized
by DAEs, such as the approaches in \cite{morGosZA20,morAntGH20}.

In this class, the linearized Navier-Stokes equations (NSEs) are a meaningful
test case. The standard velocity-pressure formulation of the NSEs comes with a
matrix-pencil of index 2; however, the linear term in the transfer function only
occurs if the control happens to affect the continuity equation; cp.
\cite{morAhmBGetal17}. In a theoretical model, this scenario hardly ever occurs
from first principles. Nonetheless, in the numerical realization of Dirichlet
conditions, such a control term may appear; cp. \cite{benner2015time}. In this
work, we present a discretization approach that explicitly treats the Dirichlet
control such that the linear part in the transfer function becomes a significant
part of the model. We test established and recent developments of the LF on this
model in terms of the qualitative identification of the non-proper terms and the
quantitative approximation of the overall transfer function. 

As the title suggests, the explicit approach for matching of non-proper transfer functions in the LF will be made according to \cite{morAntGH20}. There, the polynomial coefficients are explicitly estimated and then, the data are post-processed so that a proper rational function can be fitted instead. Then, the implicit approach consists of an adaptation of the AA method in \cite{AA86} to account for such polynomial terms, directly in the barycentric form of the fitted rational function. The latter has the advantage that no estimation of polynomial coefficient or post-processing is required. 

\section{The tangential rational interpolation problem}

For MIMO (multi-input multi-output) systems, the samples of the transfer function $\bH(s) = \bC(s\bE-\bA)^{-1}\bB+\bD$ are $p\times m$ matrices. So, in the case of rational matrix interpolation, one possibility is to interpolate along specific directions otherwise, the dimension will scale with the lengths of the input-output spaces. To avoid this, a viable way is to address the so-called \textit{tangential interpolation problem} (see, e.g., \cite{antoulas2017tutorial,karachalios2021loewner}). 

We are given a set of \textit{input/output} response measurements characterized by left interpolation points $\{\mu_{i}\}_{i=1}^{q}\subset\IC$, using left tangential directions $\{\bell_{i}\}_{i=1}^{q}\subset\IC^{p}$, and producing left responses $\{\bv_{i}\}_{i=1}^{q}\subset\IC^{m}$, together with right interpolation points $\{\lambda_{i}\}_{i=1}^{k}\subset\IC$, using right tangential directions: $\{\br_{i}\}_{i=1}^{k}\subset\IC^{m}$, producing right responses: $\{\bw_{i}\}_{i=1}^{k}$.

We are thus given the left data subset $(\mu_{j};{\bell_{j}}^{T}, \bv_j^{T})$, $j=1,\ldots ,q$, and also 
the right data subset $(\lambda_{i};\br_{i}, \bw_i)$, $i=1,\ldots ,k$. The goal is to find a rational $p\times m$ matrix function $\bH(s)$, such that the tangential interpolation conditions below are matched:
\begin{align} \label{TangIP}
	\begin{split}
	\bH(\lambda_i)\br_i&=\bw_i,~i=1,\ldots,k,\\
	\bell_j^T\bH(\mu_j)&=\bv_j^T,~
	j=1,\ldots,q.
	\end{split}
\end{align}
The \underline{left data subset} is interpolation points rearranged as:
\begin{equation}\label{leftdata}
	\bM=
\text{diag}(\mu_1, \ldots,\mu_q) \in\IC^{q\times q},
	\begin{array}{l}
		\bL^T=[\bell_1~~\cdots~~\bell_q] \in\IC^{p\times q}, \\[2mm]
		\IV^T=[\bv_1~~\cdots~~\bv_1]\in\IC^{m\times q}.
	\end{array}
\end{equation}
while the \underline{right data subset} is arranged as:
\begin{equation}\label{rightdata}
	\bLambda=
		\text{diag}(\lambda_1, \ldots,\lambda_k) \in\IC^{k\times k},
	\begin{array}{l}
		\bR=[\br_1~~\cdots~~\br_k] \in\IC^{m\times k}, \\[2mm]
		\IW=[\bw_1~~\cdots~~\bw_k]\in\IC^{p\times k}.
	\end{array}
\end{equation}
Interpolation points and tangential directions are determined by the problem or are selected to realize given MOR goals. 
\begin{equation}\label{mimocond}
	\left.\begin{array}{l}
		\bell_j^T\hat{\bH}(\mu_{j})=\bell_j^T\bH(\mu_{j})\Rightarrow \bell_j^T\hat{\bH}(\mu_{j})=\bv_{j},~
		j=1,\cdots,q,\\[1mm]
		\hat{\bH}(\lambda_{i})\br_i=\bH(\lambda_{i})\br_i\Rightarrow \hat{\bH}(\lambda_{i})\br_i=\bw_{i},
		~i=1,\cdots,k.
	\end{array}\right\}
\end{equation}

For SISO systems, i.e., $m=p=1$, left and right directions can be taken equal to one ($\bell_{j}=1, \br_{i}=1$), and hence the conditions above become:
\begin{equation}\label{sisocond}
	\left.\begin{array}{l}
		\hat{\bH}(\mu_{j})=\bH(\mu_{j})\Rightarrow \hat{\bH}(\mu_{j})=\bv_{j},~
		j=1,\cdots,q,\\[1mm]
		\hat{\bH}(\lambda_{i})=\bH(\lambda_{i})\Rightarrow \hat{\bH}(\lambda_{i})=\bw_{i},
		~i=1,\cdots,k.
	\end{array}\right\}
\end{equation}

\subsection{Interpolatory projectors}
For arbitrary values $k,q$, the matrices $\cR 	\in\IC^{n\times k}$ and $\cO^T	\in\IC^{n\times k}$ defined below
\begin{equation} \label{proj11}
	\cR=\left[(\lambda_1\bE-\bA)^{-1}\bB\br_1,~\cdots~,~
	(\lambda_k\bE-\bA)^{-1}\bB\br_k\right],
\end{equation}
\begin{equation} \label{proj33}
	\cO^T=\left[(\mu_1\bE^T -\bA^T)^{-1}\bC^T\ell_1~~\cdots~~
	(\mu_q\bE^T -\bA^T)^{-1}\bC^T\ell_q\right]
\end{equation}
will be used as \textit{projection matrices}, in order to impose the tangential interpolation properties introduced above.

The projected system is computed via a double-sided projection-based approach, as below:
\begin{align}
	\hat\bE&=\cO\bE\cR \in\IC^{q\times k}, \ \ \hat\bA= \cO\bA\cR \in\IC^{q\times k}, \\ \hat\bB&=\cO\bB \in\IC^{k\times m}, \ \ 	\hat\bC=\bC\cR \in\IC^{p\times q}.
\end{align}

\section{The Loewner framework}

It holds that, by following the derivations in \cite{mayo2007framework}, the reduced quantities $\hat\bE$ and $\hat\bA$ form a \textit{Loewner pencil}: 
\begin{equation}
	\hat\bE=
	-\left[\begin{array}{ccc}
		\frac{\bv_1^T\br_1-\ell_1^T\bw_1}{\mu_1-\lambda_1} & \cdots &
		\frac{\bv_1^T\br_k-\ell_1^T\bw_k}{\mu_1-\lambda_k} \\
		\vdots & \ddots & \vdots \\
		\frac{\bv_q^T\br_1-\ell_q^T\bw_1}{\mu_q-\lambda_1} & \cdots &
		\frac{\bv_q^T\br_k-\ell_q^T\bw_k}{\mu_q-\lambda_k} \\
	\end{array}\right] 
	:=-\IL,
\end{equation}
\begin{equation*}
	\hat\bA=-\left[\begin{array}{ccc}
		\frac{{\mu_1}\bv_1^T\br_1-\ell_1^T\bw_1\lambda_1}{\mu_1-\lambda_1} & \cdots &
		\frac{{\mu_1}\bv_1^T\br_k-\ell_1^T\bw_k\lambda_k}{\mu_1-\lambda_k} \\
		\vdots & \ddots & \vdots \\
		\frac{\mu_q\bv_q^T\br_1-\ell_q^T\bw_1\lambda_1}{\mu_q-\lambda_1} & \cdots &
		\frac{\mu_q\bv_q^T\br_k-\ell_q^T\bw_k\lambda_k}{\mu_q-\lambda_k} \\
	\end{array}\right]
	:=-\sIL,
\end{equation*}
and also that
\begin{align*}
	\hat\bB&=\left[\begin{array}{c}
		\bv_1^T\\\vdots\\\bv_q^T\end{array}\right]:=\IV, \ \ 
	\hat\bC=\left[\begin{array}{ccc}
		\bw_1&~\cdots~&\bw_k\end{array}
	\right]:=\IW.
\end{align*}
The resulting collection of data matrices ~$(\IW,\,\IL,\,\sIL,\,\IV)$~ 
is known as the \textit{Loewner quadruple}.

\begin{lemma}
The following relations hold true:
	\begin{equation}\label{projeq}
		\sIL-\IL\,\bLambda=\IV\,\bR\,~~~\text{and}~~~
		\,\sIL-\bM\,\IL=\bL\,\IW .
	\end{equation}
Then, it directly follows that the 
	Loewner quadruple satisfies the Sylvester equations
	\begin{equation*}
		\bM\IL-\IL\,\bLambda=\IV\bR-\bL\IW, ~
		\bM\,\sIL-\sIL\bLambda=\bM\IV\bR-\bL\IW\bLambda.
	\end{equation*}
\end{lemma}

\begin{theorem}
	Assume that $q = k$ and that the pencil $(\sIL,~\IL)$ is regular\footnote{The pencil $(\sIL,\IL)$ is regular if there exists $\zeta \in \IC$  such that $ \det(\sIL-\zeta \IL)\neq 0$.}. Then $\bH(s)=\IW(\sIL-s\IL)^{-1}\IV,$
	satisfies the tangential interpolation condition (\ref{TangIP}).
\end{theorem}

\textbf{Proof:} see \cite{antoulas2017tutorial} for the precise arguments.

Parameterization of all interpolants can be achieved by artificially including a term $\bK$ as shown in the next result.
\begin{remark}
	The Sylvester equation for $\IL$ can be rewritten as
	$
	\begin{array}{c}
		\bM\IL-\IL\bLambda=(\IV-\bL\bK)\bR-\bL(\IW-\bK\bR),
	\end{array} 
	$ where $\bK\in\IC^{p\times m}$ together with a similar one for $\bar{\IL}_s$. Hence, 
 $\left(\bar{\IW},\IL,\bar{\ \IL}_s,\bar{\IV}\right)$ is an interpolant for all $\bK\in\IC^{p\times m}$, where
	$$
	\bar{\ \IL}_s=\sIL+\bL\bK\bR, ~\bar\IV=\IV-\bL\bK, ~
	\bar\IW=\IW-\bK\bR.
	$$
\end{remark}

\subsection{Construction of interpolants}

If the pencil ~$(\sIL,\,\IL)$~ is regular, then ~
$\bE=-\IL,~~ \bA=-\sIL,~~ \bB=\IV,~~ \bC=\IW$,~ 
is a minimal interpolant of the data, i.e., 
$\bH(s)=\IW(\sIL-s\IL)^{-1}\IV$, interpolates the data, and it is of minimal degree.

Otherwise, problem (\ref{TangIP}) 
has a solution \cite{mayo2007framework} provided that
\begin{equation*}\label{assumption}
	\hspace*{-2mm}
	\mbox{rank}\,\left[s\,\IL-\sIL\right]=\mbox{rank}\,\left[\IL,\ \, \sIL\right]=
	\mbox{rank}\,\left[\!\begin{array}{c}\IL\\\sIL\end{array}\!\right]\!= {r},
\end{equation*}
for all ~$s\in\{\lambda_j\}\cup\{\mu_i\}$.~ 
Consider, then, the short SVDs: 
\begin{equation*}\label{prmat}
	\left[\IL,\ \, \sIL\right]=\bY\widehat{\Sigma}_{ {r}}\tilde{\bX}^*,~~
	\left[\begin{array}{c}\IL\\\sIL\end{array}\right] = {\tilde\bY}\Sigma_{ {r}} \bX^*,
\end{equation*}
where ~$\widehat{\Sigma}_{ {r}}$, $\Sigma_{ {r}}$ $\in$ $\IR^{{ {r}}\times  {r}}$,~
$\bY \in\IC^{q\times  {r}}$,
$\bX$ $\in$ $\IC^{k \times  {r}}$, $\tilde{\bY} \in\IC^{2q\times  {r}}$, $\tilde{\bX}$ $\in$ $\IC^{r \times  {2k}}$.  

\begin{remark}
In practical applications,	the value {$r$} can be chosen as the {numerical rank} of the Loewner pencil, based on a tolerance value $\tau >0$.
\end{remark} 

\begin{theorem}
	The quadruple ~$(\bhE ,\bhA ,\bhB ,\bhC )$~ of size 
	~$ {r}\times  {r}$, ~$ {r}\times  {r}$, ~$ {r}\times m$, ~$p\times  {r}$,~ given by:
	\begin{equation*} \label{redundant}
		\bhE  = -\bY^T\IL \bX ,~~\bhA  = -\bY^T\sIL \bX, ~~\bhB  = \bY^T\IV,~~
		\bhC  = \IW \bX ,
	\end{equation*}
	is a descriptor realization of an (approximate) interpolant of the data with McMillan degree
	$~ {r}=\mbox{rank}\,\IL$. 
\end{theorem}
The approximation relies on the choice of the value $\tau$; if this value is
indeed $0$ (exact arithmetic), then the interpolation is also exact. However,
for $\tau = 0$, the interpolation errors are proportional to $\tau$ (or to the
first neglected singular value of the Loewner matrix), as shown in \cite{antoulas2017tutorial}. 
\begin{remark}
	(a) The Loewner framework constructs a 
	descriptor representation $(\IW,\IL,\sIL,\IV)$ from data, with no further processing. However, if the pencil $(\sIL,\IL)$ is singular, it needs to be projected to a regular pencil $(\bA,\bE )$. \\[1mm]
	(b) By construction, the $\bD$ term is absorbed into the other matrices of
	the realization. Extracting the $\bD$ term 
	involves an eigenvalue decomposition of $(\sIL,\IL)$, and a careful numerical treatment of the spectrum at $\infty$. \\[1mm]
	(c) If the transfer function of the underlying system has higher-order polynomial terms, then a more intricate procedure is needed to accurately recover these terms (behavior at infinity) in the LF; some solutions were recently proposed in \cite{morGosZA20,morAntGH20}.
\end{remark}

In the sequel, we follow the approach originally proposed in \cite{morAntGH20}, for estimating the polynomial terms from values of the transfer function, sampled at so-called high frequencies. There, what is proposed is to subtract the fitted polynomial part from the original data, and perform classical Loewner framework analysis (as in \cite{mayo2007framework}) on the pre-processed data.

\subsection{Estimating the polynomial terms from data in the LF
\cite{morAntGH20}}\label{sec:poly-loewner}

We assume that the transfer function of the underlying large-scale model is
composed of a strictly proper part and of a polynomial part, in the following
way $\bH(s) =  \bHspr(s) + \bHpoly(s)$. The polynomial part is considered to have only two non-zero coefficients (and as such, to be a linear polynomial in $s$). In the widely-accepted terminology, this scenario corresponds to a DAE of index 2, with the polynomial part as below:
\begin{equation}
	\bHpoly(s) = \bP_0 + s \bP_1, \ \ \ \bP_0, \ \bP_1 \in \IR^{p \times m}.
\end{equation}
The main idea summarized in \cite{morAntGH20}, is that by having access at
$\bH(s)$ for large values of $s$, the contribution of $\bHspr(s)$ to the transfer function is negligible; hence, it can be ignored. We review some of the formulae presented in the aforementioned contribution, first for only limited amount of data, and afterwards, for many data points. The estimates of the two coefficients will be denoted with $\bhP_i$, for $\bP_i$, where $0 \leq i \leq 1$. Also, for the first cases, no tangential directions will be used.

\subsection{The case of a few data points}

We assume that the transfer function is known at one value, $\imath \omega_1$, where $\imath := \sqrt{-1}$ and $\omega \in \IR$. Then, it holds that:
\begin{equation}\label{P0_P1_approx-1value}
	\bhP_0 = \Re(\bH(\imath \omega_1), \ \
	\bhP_1 = \omega^{-1} \Im(\bH(\imath \omega ).
\end{equation}
Then, if two sample values are known, i.e., at the points $j \omega_1$ and $j \omega_2$ on the imaginary axis, with $\omega_k \in \IR$ for $1 \leq k \leq 2$. Then, as shown in \cite{morAntGH20}, the following estimates hold true:
\begin{align}\label{P0_P1_approx-2value}
	\bhP_0 &= \Re \Big{(}\frac{\imath \omega_1 \bH(\imath \omega_1)-\imath \omega_2 \bH(\imath \omega_2)}{\imath \omega_1- \imath \omega_2}\Big{)}, \\
	\bhP_1 &= \frac{\bH(\imath \omega_1)-\bH(\imath \omega_2)}{\imath \omega_1- \imath \omega_2}.
\end{align}

\subsection{The general case (many data points)}

Assume now that $k=q \geq 2$ and that $2k$ interpolation points, values, and tangential directions are provided.
Instead of the generic notation $\IL$ for the  Loewner matrix, we now use the
notation $\IL^{\mathsf{hi}}$ to 
indicate that this  Loewner matrix is computed with
data located in high frequency ranges.

Provided that $k \geq \max\{p,m\}$, one can write the estimated linear polynomial coefficient matrix as shown in \cite{morAntGH20}, as:
	\begin{equation}\label{P1_approx_gen}
    \bhP_1 = \big(\bL^{\mathsf{ hi}}\big)^\dagger \IL^{\mathsf{ hi}}
    \big(\bR^{\mathsf{ hi}} \big)^\dagger,
	\end{equation}
	where $\bX^\dagger \in \mathbb{C}^{v \times u}$ is the Moore-Penrose pseudo-inverse of the
	matrix $\bX \in \mathbb{C}^{u \times v}$.
	
	Similarly to the procedure used for estimating $ \bP_1$, one can extend the formula  for estimating $ \bP_0$ from the shifted Loewner matrix 
	$ \IL_s^{\mathsf{hi}}$ computed from $2k$ sampling points  located in high frequency bands
	as follows (as shown in \cite{morAntGH20})
	\begin{equation}\label{P0_approx_gen}
		\bhP_0 = \Re\Big( \big(\bL^{\mathsf{hi}}\big)^\dagger \IL_s^{\mathsf{hi}}  \big(\bR^{\mathsf{hi}} \big)^\dagger \Big).
	\end{equation}
	
\section{Other methods}	

Here, we will only go into details for direct methods, which do not require an iteration (for a fair comparison with \cite{morAntGH20}). It is to be mentioned that the Vector Fitting algorithm in \cite{morGusS99} is a robust, effective rational approximation tool based on a least-square fit on the data, and also can be used to accommodate up to linear polynomial terms (of the fitted transfer function). However, we will not concentrate on this here.

\subsection{Antoulas-Anderson method with (higher) polynomial
terms}\label{sec:poly-AA}

The rational approximant $\hat{\bH}(s)$ computed by the original Antoulas-Anderson (AA) rational approximation approach in \cite{AA86} is based upon the classical barycentric form:
\begin{equation}\label{eq:Barymodel}
	\hat{\bH}(s)=  \frac{\sum_{i=1}^k \displaystyle \frac{w_i h_i}{s-z_i}}{\sum_{i=1}^k \displaystyle \frac{w_i}{s-z_i}}=\frac{N(s)}{D(s)}.
\end{equation}
This representation generally enforces proper ($\text{deg}(N(s))\leq \text{deg}(D(s))$) or strictly proper ($\text{deg}(N(s)) < \text{deg}(D(s))$) transfer functions. It is to be noted that the constant polynomial term can be recovered as the ratio of sums $\bhP_0 = \frac{\sum_{i=1}^k w_i h_i}{\sum_{i=1}^k w_i}$.

However, we are interested in recovering improper transfer functions, that can be written as proper rational functions with linear polynomial terms. Hence, the following representation of the fitted transfer function will be used instead:
\begin{equation}\label{eq:Barymodel2}
	\hat{\bH}_{\mathsf{mod}}(s)=  \frac{b+\sum_{i=1}^k \displaystyle \frac{w_i h_i}{s-z_i}}{\sum_{i=1}^k \displaystyle \frac{w_i}{s-z_i}}.
\end{equation}
Clearly, the following interpolation conditions are enforced, solely by the \textit{barycentric structure}, as in the usual case (\ref{eq:Barymodel}):
\begin{equation}
	\hat{\bH}_{\mathsf{mod}}(z_i) = h_i, \ \forall \ 1 \leq i \leq k.
\end{equation}
In this context, the variables to be fitted are the weights $w_1, w_2, \ldots, w_k$, but also the free term $b$ in the numerator.

As in in \cite{AA86}, to enforce additional $q$ (left) interpolation conditions given by:
\begin{equation}
		\hat{\bH}_{\mathsf{mod}}(s_j) = g_j, \ \forall \ 1 \leq j \leq q,
\end{equation}
one can write the problem explicitly for any $1\leq j \leq q$, as:
\begin{align}
		\hat{\bH}_{\mathsf{mod}}(s_j)&= g_j\Leftrightarrow\frac{b+\sum_{i=1}^{k}\frac{w_i h_i}{s_j-z_i}}{\sum_{i=1}^{k}\frac{w_i}{s_j-z_i}}=g_j,~\forall 1\leq j\leq q \\
	&\Leftrightarrow \sum_{i=1}^{k}\frac{g_j-h_i}{s_j-z_i}w_i-b=0,
\end{align}
or equivalently, in matrix format, as:
\begin{equation}\label{addit_interp}
\left[\begin{array}{ccccc}
		\frac{g_1-h_1}{s_1-z_1} & 	\frac{g_1-h_2}{s_1-z_2} & \cdots & 	\frac{g_1-h_k}{s_1-z_k} & -1\\[1mm]
		\frac{g_2-h_1}{s_2-z_1} & 	\frac{g_2-h_2}{s_2-z_2} & \cdots & 	\frac{g_2-h_k}{s_2-z_k} & -1\\[1mm] 
		\vdots & \vdots & \ddots & \vdots & \vdots \\[1mm]
		\frac{g_q-h_1}{s_q-z_1} & 	\frac{g_q-h_2}{s_q-z_2} & \cdots & 	\frac{g_q-h_k}{s_q-z_k} & -1 \\
	\end{array}\right] \begin{bmatrix}
	w_1 \\ w_2 \\ \vdots \\ w_k \\ b
	\end{bmatrix} = \begin{bmatrix}
	0 \\ 0 \\ \vdots \\ 0 \\ 0
	\end{bmatrix}.
\end{equation}
This can further be simply written as $\tilde{\ \IL} \bfa = \bfz$, where $\tilde{\ \IL}$ is the Loewner matrix with a column of $-1$'s augmented at the end as in (\ref{addit_interp}) and $\bfa$ is the vector of variables, i.e., $\bfa = \begin{bmatrix}
	w_1 & w_2 & \cdots & w_k & b
\end{bmatrix}^T$.
The $k+1$ unknowns in vector $\bfa$ can be computed from the (approximate) null
space, i.e., the kernel of the  matrix $\tilde{\ \IL}$.  It is to be noted that
the SVD of matrix $\tilde{\ \IL}$ can be employed for this task. 

In the case of not enough data, or noisy/perturbed measurements, the last singular value of $\tilde{\ \IL}$ is seldom zero. By setting up a tolerance value $\tau$, one can compute $\bfa$ as the left singular vector of $\tilde{\ \IL}$ corresponding to the smallest singular value greater than $\tau$. By finding the missing coefficients in vector $\bfa$, the fitted rational function will be uniquely determined. Afterwards, if needed, the polynomial terms of $\hat{\bH}_{\mathsf{mod}}(s)$ can be explicitly explicitly in terms of the recovered coefficients in vector $\bfa$, the interpolation points, and transfer function measurements. We skip this step here for brevity reasons.

\section{Numerical Examples}

\subsection{A first example}

 We use the damped mass-spring system with a holonomic constraint example (in
 short, \MSD) from \cite{morMehS05}. Although originally an index-3 DAE system,
 it can be transformed into an index-1 DAE system by appropriately choosing the
 $\bB$ and $\bC$ vectors. We denote the resulting full model transfer function
 by $H_{\mathsf{MSD}}$ and refer the reader to the original publication for more details.

\subsection{Oseen equations with Dirichlet control}

We consider a flow example with boundary control modeled by a finite element
discretization of the incompressible Oseen equations. 
The Oseen equations are obtained from the Navier-Stokes equations by a Newton
linearization about a steady state solution.
We will consider setups that, apart from the boundary $\Gamma$ where the control acts,
have \emph{no-slip} boundary conditions at the walls or \emph{do-nothing} conditions at the outlets. 

The control $\nu(t, x)$, where $t$ denotes the time and $x$ denotes the
spatial variable distributed over the considered boundary, is modeled as
$\nu(t, x) = g(x)u(t)$ through a
function $g\colon \Gamma \to \mathbb R^{2/3}$ that describes the spatial
extension and through a scalar function $u$ that models the control action as a
scaling of $g$.

Overall, the spatially-discretized model for the velocity $v$ and
pressure $p$ reads
\begin{equation}\label{eq:oseen-semidisc-system}
  \begin{split}
    \begin{bmatrix}
      M & M_\Gamma
    \end{bmatrix}
    \begin{bmatrix}
      \dot v(t) \\ \dot v_\Gamma(t)
    \end{bmatrix}
    &=
    \begin{bmatrix}
      A & A_\Gamma
    \end{bmatrix}
    \begin{bmatrix}
      v(t) \\ v_\Gamma(t)
    \end{bmatrix}
    + J\tt  p(t) , \\
    0 &= 
    \begin{bmatrix}
      J & J_\Gamma
    \end{bmatrix}
    \begin{bmatrix}
      v(t) \\ v_\Gamma(t)
    \end{bmatrix},
    \\
    0 & =v_\Gamma (t) - b_\Gamma u(t),
  \end{split} 
\end{equation}
where $v_\Gamma$ denotes the degrees of freedom in the discrete velocity vector
that are associated with the control boundary, where the corresponding parts
of the linear operators are subscripted with $\Gamma$ accordingly, and where
$b_\Gamma$ is the spatially discretized representation of $g$.

If one resolves $v_\Gamma(t)=b_\Gamma u(t)$ directly, the following controlled linear
system is obtained
\small
\begin{subequations}\label{eq:bc-oseen-resolved}
  \begin{align}
  M\dot v(t) &= Av(t) + J\tt p(t) + f(t) + A_\Gamma b_\Gamma u(t) -
  M_\Gamma b_\Gamma  \dot u(t) \\
  0 &= Jv(t) - g(t) - J_\Gamma b_\Gamma u(t)
\end{align}
\end{subequations}
\normalsize
which we will write as
\begin{equation}
  \mathcal  E\dot x(t) = \mathcal A x(t) +\mathcal  B_1u(t) +\mathcal  B_2\dot u (t),
\end{equation}
with $x=(v,p)$ and
\begin{equation}
  \mathcal E:=
\begin{bmatrix}
  M & 0 \\ 0 & 0
\end{bmatrix}, \quad
  \mathcal A:=
\begin{bmatrix}
  A & J\tt \\ J & 0
\end{bmatrix},
\end{equation}
and
\begin{equation}
\mathcal B_1 =
\begin{bmatrix}
  A_\Gamma b_\Gamma \\ J_\Gamma b_\Gamma
\end{bmatrix}, \quad
\mathcal B_2 =
\begin{bmatrix}
  -M_\Gamma b_\Gamma \\ 0.
\end{bmatrix}
\end{equation}

As laid out in \cite[Sec. 5]{morAhmBGetal17}, depending on how the output 
\begin{equation}
  y(t) = C_v v(t) + C_p p(t)
\end{equation}
is defined, the transfer function $u\mapsto y$ associated with the system
\eqref{eq:bc-oseen-resolved} can be strictly proper, proper, or include a
linear part. 
In particular, if $C_p \neq 0$, i.e., the measurements include the $p$-variable,
then the transfer function will likely have a linear part.

\subsection{The particular setup}

\begin{figure}
  \includegraphics[width=\linewidth]{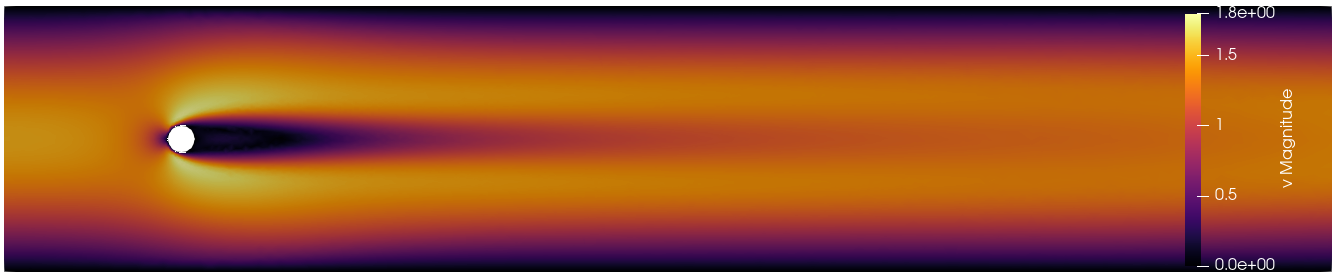}
  \caption{Snapshot of magnitude of the steady-state Navier-Stokes velocity solution in the considered setup.}
  \label{fig:stst-nse}
\end{figure}

We consider the flow past a cylinder in two dimensions at Reynolds number 20
calculated with the averaged inflow velocity and the cylinder diameter as
reference quantities; see Fig. \ref{fig:stst-nse} for a velocity magnitude snapshot.

As the control input, we consider a modulation the input velocity around its
reference value. As the observation, we consider a single output $y_p$ of the pressure
or the single output $y_v$ of the sum of the two velocity
components averaged over a square domain of area $d^2$ located $10d$ behind
the cylinder in the wake, where $d$ is the cylinder diameter.

The spatial discretization is obtained by \emph{Taylor-Hood} piecewise
quadratic/piecewise linear finite elements on an unstructured triangulation of
the domain which results in around 42000 degrees of freedom for the velocity and
5000 degrees of freedom for the pressure.

\begin{figure}
	\centering
  \includegraphics[width=0.9\linewidth]{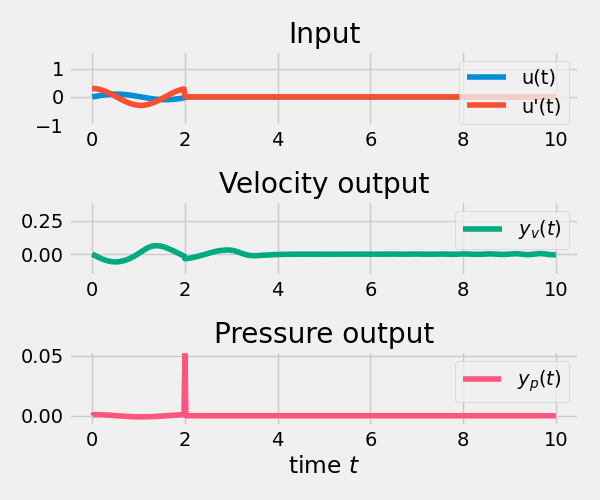}
  \caption{Time-domain linear response of the Oseen system for a test output set
  to zero for $t\geq 2$. Notably, and in line with the theory, the resulting impulse is observed in the
pressure output $y_p$ but not in the velocity output $y_c$.}
\label{fig:nse-time-dom-response}
\end{figure}

The linearized model is obtained from linearizing the system about the corresponding steady-state solution.
Thus, the obtained input/output map models the linear response to a changing
input velocity in the measurements in the wake; see Fig.
\ref{fig:nse-time-dom-response} for an illustration of the response in time domain.
 For the presented results on transfer function interpolation, we consider the
 $y_p$ output only, i.e. 
 \begin{equation}
   H_{\mathsf{OS}}(s) := \mathcal C_p (s\mathcal E - \mathcal A)^{-1}(\mathcal B_1 + s\mathcal
   B_2),
 \end{equation}
 with $\mathcal C_p = \begin{bmatrix}
   0 & C_p
 \end{bmatrix}$.

\subsection{Numerical Results}

For further reference, we denote the \emph{Antoulas-Anderson} approach with
polynomial terms (see Sec. \ref{sec:poly-AA}) by \CAA and the Loewner with
explicit matching by (see Sec. \ref{sec:poly-loewner}) by \pLWNR.

We investigate the performance of the \CAA~approach for the \MSD~and the
\Oseen~example by comparing to a plain Loewner interpolation and the Loewner
approach with explicit identification of the constant and the linear term in Sec. \ref{sec:poly-loewner}. The parameters for the approximations are chosen as follows:
\begin{itemize}
  \item For determining the polynomial part in \pLWNR, we consider the range
    $[10^7j,\, 10^9j]\subset j\mathbb R$ and 20 or 10 for the \MSD~or
    \Oseen~example, respectively, evenly (on the logarithm-scale) distributed interpolation points.
  \item For determining the proper part in the plain Loewner (and the \pLWNR),
    we consider 200 or 40 interpolation points (for \MSD~or \Oseen~respectively)
    evenly log-distributed on $[10^{-2}j,\, 10^4j]\subset j\mathbb R$.
  \item The order $r$ of the identified is inferred by truncating all singular
    values smaller than the relative tolerance $10^{-10}$.
  \item To compute the \CAA~interpolation, for both setups, we used 120
    interpolation points, defined as 48 and 12 \emph{left} and \emph{right}
    interpolation points evenly distributed on $[10^{-3}j,\,10^6j]\subset
    j\mathbb R$ plus their complex conjugates.
\end{itemize}

\begin{figure}
  \includegraphics[width=\linewidth]{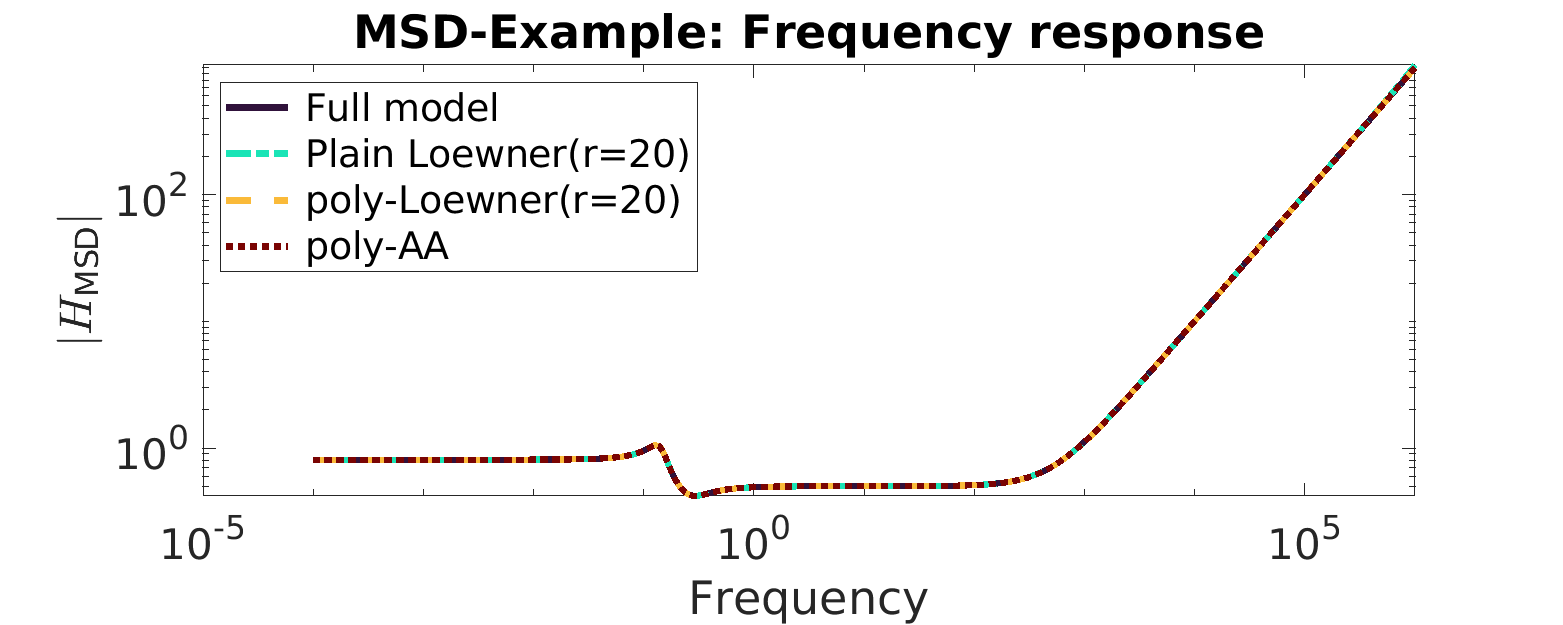}
  \includegraphics[width=\linewidth]{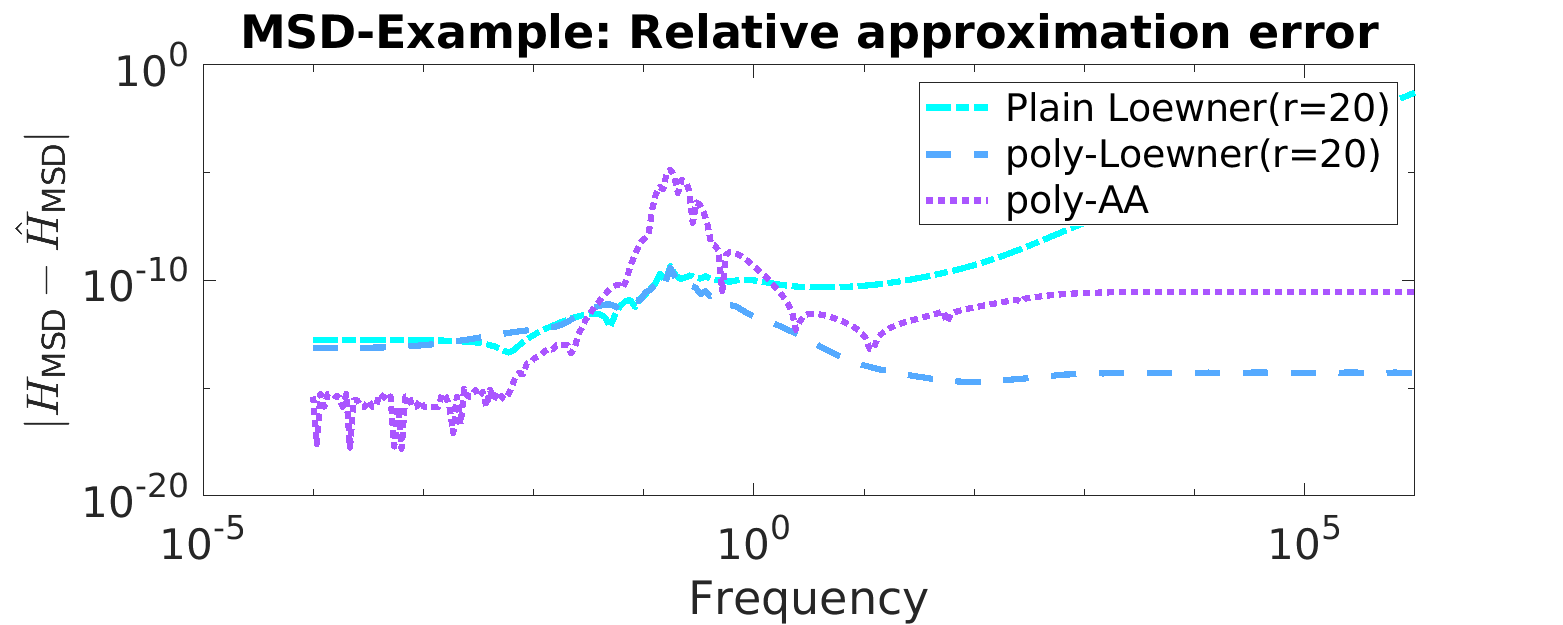}
  \caption{Frequency responses of the \MSD~full order model and the
  interpolations (top) and the relative errors (bottom).}
  \label{fig:freqrep-relerr-msd}
\end{figure}

\begin{figure}
  \includegraphics[width=\linewidth]{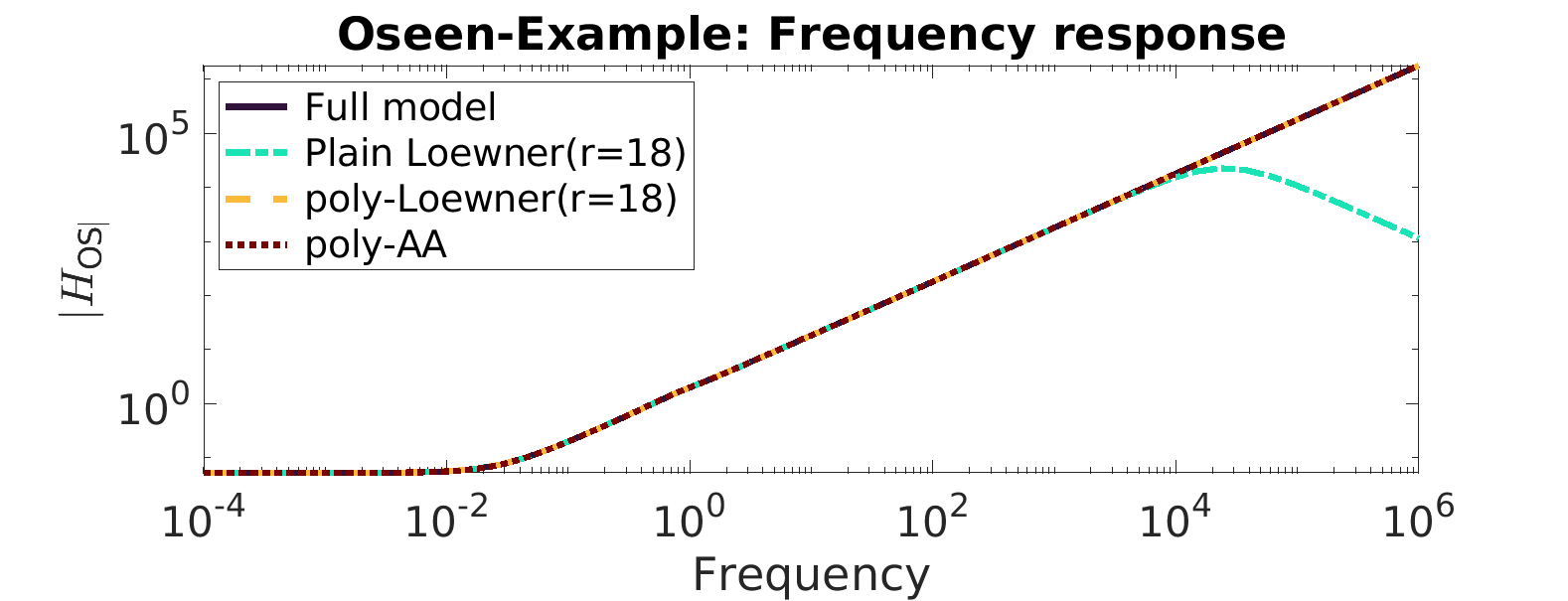}
  \includegraphics[width=\linewidth]{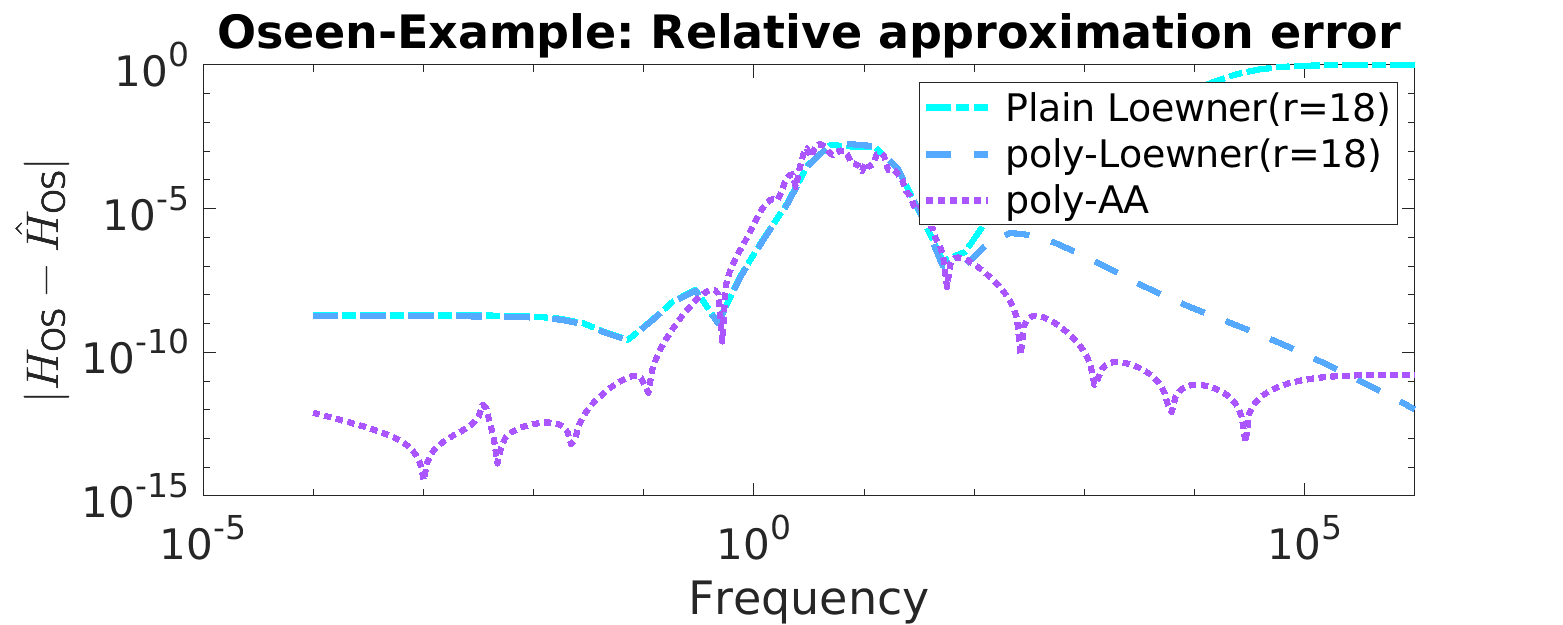}
  \caption{Frequency responses of the \Oseen~full order model and the
  interpolations (top) and the relative errors (bottom).}
  \label{fig:freqrep-relerr-oseen}
\end{figure}

The simulation results are plotted in Fig. \ref{fig:freqrep-relerr-msd}
(\MSD~example) and Fig. \ref{fig:freqrep-relerr-oseen} (\Oseen~example).

In both examples, the frequency response of the full model and approximations
are visually indistinguishable. 
The plots of the relative errors,, however, reveal the qualitative
difference between the plain Loewner approach (which does not capture the linear
behavior that dominates for high frequencies) and \pLWNR~and~\CAA{} as well as
quantitative differences between \pLWNR{} and~\CAA. 
While the good performance of
\CAA~in the low-frequency regime (and the indifferent performance in the
mid-frequency range) are likely be explained by the number and choice of
interpolation points, we attribute the reliably worse performance of \CAA~for
high-frequencies to the all-at-once determination of proper and non-proper
components (cp. Eq. \eqref{addit_interp}) of the transfer function. 


\subsection*{Code Availability}
The linear system data (in the form of a \texttt{.mat} file) and the scripts that were
used to obtain the presented numerical results are available for immediate
reproduction from 
\href{https://dx.doi.org/10.5281/zenodo.10058537}{doi:10.5281/zenodo.10058537}
under a CC-BY license.

\section{Conclusion}

We have proposed a variant of Loewner-based system identification based on the
\emph{Antoulas-Anderson} algorithm but with free parameters to account for
non-proper parts. 
In this approach, polynomial parts of the transfer function are implicitly covered,
which is advantageous over explicit treatments of the polynomial parts that requires data points at sufficiently large frequencies.
As a drawback of the implicit realization, the determination through, say,
singular value decompositions, gives little error control on the individual
coefficients and, thus, lead to a larger approximation error in the
high frequency range. On the other hand, the equal treatment of all
interpolation points gives way to consider adaptive versions of the
\CAA~approach as in the \emph{adaptive Antoulas-Anderson} algorithm; see \cite{morNakST18}.
Another future investigations will concern extensions of the \CAA~approach to possibly
higher polynomial terms.

\bibliographystyle{IEEEtran}
\bibliography{nonproper-nse-loewner}

\end{document}

%% file: my_macros.tex
\def\IR{{\mathbb R}}
\def\IC{{\mathbb C}}

\def\IL{{\mathbb L}}

\def\IV{{\mathbb V}}
\def\IW{{\mathbb W}}

\newcommand{\sIL}{{{{\mathbb L}_s}}}

\newcommand{\bA}{{\mathbf A}}
\newcommand{\bB}{{\mathbf B}}
\newcommand{\bC}{{\mathbf C}}
\newcommand{\bD}{{\mathbf D}}
\newcommand{\bE}{{\mathbf E}}

\newcommand{\bY}{{\mathbf Y}}

\newcommand{\bK}{{\mathbf K}}

\newcommand{\bL}{{\mathbf L}}
\newcommand{\bM}{{\mathbf M}}

\newcommand{\bH}{{\mathbf H}}
\newcommand{\bP}{{\mathbf P}}

\newcommand{\bR}{{\mathbf R}}

\newcommand{\bX}{{\mathbf X}}

\newcommand{\bfa}{{\mathbf a}}

\newcommand{\br}{{\mathbf r}}

\newcommand{\bv}{{\mathbf v}}
\newcommand{\bw}{{\mathbf w}}

\newcommand{\bfz}{{\mathbf 0}}

\newcommand{\cO}{ {\mathcal O} }

\newcommand{\bell}{\boldsymbol{\ell}}

\newcommand{\bLambda}{\boldsymbol{\Lambda}}

\newcommand{\cR}{ {\cal R} }

\newcommand{\bhA}{{\hat{\textbf A}}}
\newcommand{\bhB}{{\hat{\textbf B}}}
\newcommand{\bhC}{{\hat{\textbf C}}}

\newcommand{\bhE}{{\hat{\textbf E}}}

\newcommand{\bhP}{{\hat{\textbf P}}}

\renewcommand{\Re}{\operatorname{Re}}
\renewcommand{\Im}{\operatorname{Im}}

\newcommand{\CAA}{\texttt{poly-AA}}
\newcommand{\MSD}{MSD}
\newcommand{\Oseen}{Oseen}
\newcommand{\pLWNR}{\texttt{poly-Loewner}}

\newtheorem{theorem}{Theorem}[section]

\newtheorem{lemma}[theorem]{Lemma}
\newtheorem{remark}[theorem]{Remark}